\documentclass{amsart}
\usepackage{amscd,amssymb,enumerate,amsmath}
\usepackage[all]{xy}

\newcommand{\ppq}{\leqslant}
\newcommand{\pgq}{\geqslant}
\newcommand{\op}{\oplus}

\newcommand{\car}{\operatorname{char}\nolimits}

\newcommand{\Ext}{\operatorname{Ext}\nolimits}
\newcommand{\HH}{\operatorname{HH}\nolimits}
\newcommand{\opp}{\operatorname{op}\nolimits}
\newcommand{\gr}{{\operatorname{gr}\nolimits}}

\newcommand{\soc}{\operatorname{soc}\nolimits}
\renewcommand{\op}{\operatorname{op}\nolimits}
\renewcommand{\sp}{\operatorname{sp}\nolimits}

\renewcommand{\lq}{\Lambda_{\mathbf q}}
\newcommand{\iq}{I_{\mathbf q}}
\newcommand{\G}{\Gamma}
\newcommand{\oa}{\bar{a}}
\newcommand{\mo}{\mathfrak{o}}
\newcommand{\mt}{\mathfrak{t}}
\newcommand{\rrad}{\mathfrak{r}}

\newcommand{\N}{\mathcal{N}}

\newcommand{\set}[1]{\left\{ #1 \right\}}

\newtheorem{theorem}{{Theorem}}[section]
\newenvironment{thm}{\begin{theorem}}{\end{theorem}}
\newcommand{\bt}{\begin{thm}}
\newcommand{\et}{\end{thm}}

\newtheorem{corollaire}[theorem]{{Corollary}}
\newenvironment{cor}{\begin{corollaire}}{\end{corollaire}}
\newcommand{\bc}{\begin{cor}}
\newcommand{\ec}{\end{cor}}

\newtheorem{lemme}[theorem]{{Lemma}}
\newenvironment{lemma}{\begin{lemme}}{\end{lemme}}
\newcommand{\bl}{\begin{lemma}}
\newcommand{\el}{\end{lemma}}

\newtheorem{proposition}[theorem]{{Proposition}}
\newcommand{\bprop}{\begin{proposition}}
\newcommand{\eprop}{\end{proposition}}

\newtheorem{definition}[theorem]{{Definition}}
\newenvironment{dfn}{\begin{definition} \rm}{\end{definition}}
\newcommand{\bd}{\begin{dfn}}
\newcommand{\ed}{\end{dfn}}

\newtheorem*{remark}{Remark}
\newcommand{\br}{\begin{remark}}
\newcommand{\er}{\end{remark}}

\begin{document}

\topmargin 0cm
\oddsidemargin 0.5cm
\evensidemargin 0.5cm
\pagestyle{plain}

\title[Hochschild cohomology of socle deformations \dots]
{Hochschild cohomology of socle deformations of a class of Koszul self-injective algebras}
\author[Snashall]{Nicole Snashall}\thanks{The first author would like to thank the LaMUSE at the University of St Etienne for the invited professorship held during the research for this paper.}
\address{Nicole Snashall\\Department of Mathematics\\
University of Leicester\\
University Road\\
Leicester, LE1 7RH\\
England}
\email{N.Snashall@mcs.le.ac.uk}
\author[Taillefer]{Rachel Taillefer}
\address{Rachel Taillefer\\Universit\'e de Lyon, F-42023 Saint-Etienne, France; LaMUSE, Laboratoire de Math\'ematiques de l'Universit\'e de Saint-Etienne Jean Monnet, F-42023, France.
Facult\'e des Sciences et Techniques\\
23 Rue Docteur Paul Michelon\\
42023 Saint-Etienne Cedex 2\\
France}
\email{rachel.taillefer@univ-st-etienne.fr}

\begin{abstract}
We consider the socle deformations arising from formal deformations of a class of Koszul self-injective special
biserial algebras which occur in the study of the Drinfeld double of the generalized Taft algebras. We show, for these deformations, that the Hochschild cohomology ring modulo nilpotence is a finitely generated commutative algebra of Krull dimension 2.
\end{abstract}

\date{\today}
\maketitle

\section*{Introduction}

Let $K$ be a field. For $m \pgq 1$, let ${\mathcal Q}$ be the
quiver with $m$ vertices, labelled $0, 1, \ldots, m-1$, and $2m$
arrows as follows:
$$
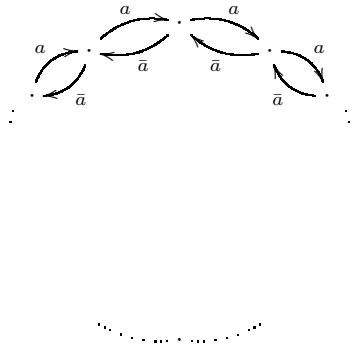
\xymatrix@=.01cm{
&&&&&&&&&&&\cdot\ar@/^.5pc/[rrrrrd]^{a}\ar@/^.5pc/[llllld]^{\bar{a}}\\
&&&&&&\cdot\ar@/^.5pc/[rrrrru]^{a}\ar@/^.5pc/[llldd]^{\bar{a}}&&&&&&&&&&\cdot\ar@/^.5pc/[rrrdd]^{a}\ar@/^.5pc/[lllllu]^{\bar{a}}\\\\
&&&\cdot\ar@/^.5pc/[rrruu]^{a}\ar@{.}@/_.3pc/[ldd]&&&&&&&&&&&&&&&&\cdot\ar@/^.5pc/[llluu]^{\bar{a}}\ar@{.}@/^.3pc/[rdd]\\\\
&&&&&&&&&&&&&&&&&&&&&&&\\
\\
\\
\\\\\\\\\\\\\\\\
&&&&&&&&&&&&&&&&\\
&&&&&&&&&&&\cdot\ar@/_.3pc/@{.}[rrrrru]\ar@/^.3pc/@{.}[lllllu] }$$
Let $a_i$ denote the arrow that goes from vertex $i$ to vertex
$i+1$, and let $\oa_i$ denote the arrow that goes from vertex $i+1$
to vertex $i$, for each $i=0, \ldots , m-1$ (with the obvious
conventions modulo $m$). We denote the trivial path at the vertex
$i$ by $e_i$. Paths are written from left to right.

We define $\Lambda$ to be the algebra $K{\mathcal Q}/I$
where $I$ is the ideal of $K{\mathcal Q}$ generated by
$a_ia_{i+1}$, $\oa_{i-1}\oa_{i-2}$ and $a_i\oa_i-\oa_{i-1}a_{i-1}$,
for $i=0, \ldots , m-1$, where the subscripts are taken modulo $m$.
These algebras are Koszul self-injective special biserial algebras and as
such play an important role in various aspects of representation
theory of algebras. In particular, for $m$ even, this algebra
occurred in the presentation by quiver and relations of the Drinfeld
double of the generalised Taft algebras studied in \cite{EGST}, and in
the study of the representation theory of $U_q(\mathfrak{sl}_2)$, for which,
see \cite{CK, Patra, Suter, Xiao}. The Hochschild cohomology ring of the algebra $\Lambda$ was determined in \cite{ST} (where $\Lambda = \Lambda_1$ in that notation).

In this paper we study socle deformations $\lq$ of the algebra $\Lambda$, where ${\mathbf q} = (q_0, q_1, \ldots , q_{m-1}) \in (K^*)^m$. The first section of this paper shows that $\lq$ arises from a formal deformation with infinitesimal in $\HH^2(\Lambda)$, and, further, that it is a socle deformation of $\Lambda$, that is, $\lq$ is self-injective and $\lq/\soc(\lq) \cong \Lambda/\soc(\Lambda)$. The algebras $\lq$ for $m=1$ were studied in
\cite{BGMS}, where they were used to answer negatively a question of Happel, in that their Hochschild cohomology ring is finite-dimensional but they are of infinite global dimension when $q\in K^*$ is not a root of unity.

For a finite-dimensional $K$-algebra $\Gamma$ with Jacobson radical $\rrad$, the Hochschild cohomology ring of $\Gamma$ is given by $\HH^*(\G) = \Ext^*_{\G^e}(\G,\G) = \oplus_{n \pgq 0}\Ext^n_{\G^e}(\G, \G)$ with the Yoneda product, where $\G^e =
\G^{\opp} \otimes_K \G$ is the enveloping algebra of $\G$. Since all tensors are over the field $K$ we write $\otimes$ for $\otimes_K$ throughout. We denote by $\N$ the ideal of $\HH^*(\G)$ which is generated by all homogeneous nilpotent elements. Thus $\HH^*(\G)/\N$ is a commutative $K$-algebra. The Ext algebra $E(\G)$ is defined by $E(\G) = \Ext^*_{\G}(\G/\rrad,\G/\rrad) = \oplus_{n \pgq 0}\Ext^n_{\G}(\G/\rrad, \G/\rrad)$. The graded centre of $E(\G)$ is denoted $Z_\gr(E(\G))$ and is generated by all homogeneous elements $z \in \Ext^n_\G(\G/\rrad, \G/\rrad)$ for which $zg=(-1)^{mn}gz$ for all $g \in \Ext^m_\G(\G/\rrad, \G/\rrad)$. The natural ring homomorphism $\HH^*(\G) \to E(\G)$ has image contained in $Z_\gr(E(\G))$; it was shown in \cite{BGSS} (and see \cite{K} for a generalization), that the image is precisely $Z_\gr(E(\G))$ when $\G$ is a Koszul algebra.

Section 2 of this paper describes explicitly the structure of $Z_\gr(E(\lq))$, the graded centre of the Ext algebra of $\lq$, for all $m \pgq 1$, ${\mathbf q} \in (K^*)^m$ and in all characteristics (Theorem \ref{thm:structure}). In the final section, we determine the Hochschild cohomology ring modulo nilpotence of the algebras $\lq$ for all ${\mathbf q} = (q_0, \ldots , q_{m-1}) \in (K^*)^m$, and show, in Theorem \ref{thm:HH}, that $\HH^*(\lq)/\N$ is a commutative finitely generated $K$-algebra of Krull dimension 2 when $q_0\cdots q_{m-1}$ is a root of unity.
It was conjectured in \cite{SS} that the Hochschild cohomology ring modulo nilpotence of any finite-dimensional algebra is always a finitely generated $K$-algebra. Although it was shown by Xu in \cite{Xu} (and see \cite{Sn}) that this conjecture is not true in general, with a counterexample being provided of a Koszul algebra that is not self-injective, the current paper gives a class of Koszul self-injective algebras where the Hochschild cohomology ring modulo nilpotence is a finitely generated $K$-algebra and the conjecture of \cite{SS} holds. This provides a further contribution to the study of the structure of the Hochschild cohomology ring for Koszul algebras.

\bigskip

\section{Socle deformations of $\Lambda$}

\bigskip

Let $\Lambda=K{\mathcal Q}/I$ where $I$ is
the ideal of $K{\mathcal Q}$ generated by $a_ia_{i+1},
\oa_{i-1}\oa_{i-2}$ and $a_i\oa_i-\oa_{i-1}a_{i-1}$, for
$i =0, 1, \ldots , m-1$. We write $\mo(\alpha)$ for the trivial path corresponding to the origin of the arrow $\alpha$, so that $\mo(a_i) = e_i$ and $\mo(\oa_i) =
e_{i+1}$. We write $\mt(\alpha)$ for the trivial path corresponding
to the terminus of the arrow $\alpha$, so that $\mt(a_i) = e_{i+1}$
and $\mt(\oa_i) = e_i$. Recall that a non-zero element $r \in
K{\mathcal Q}$ is said to be uniform if there are vertices $v, w$
such that $r = vr = rw$. We then write $v = \mo(r)$ and $w = \mt(r)$.

\bigskip

A minimal projective bimodule resolution $(P^n, \partial^n)$ for $\Lambda$ was given in \cite[Theorem 1.2]{ST}. With the notation of \cite{ST}, the projective $P^2$ is described by
$$P^2 = \bigoplus_{i=0}^{m-1}\bigoplus_{r=0}^2\Lambda\mo(g^2_{r,i})\otimes\mt(g^2_{r,i})\Lambda$$
where $$g_{0,i}^2 = a_ia_{i+1},\ g_{1,i}^2 = a_i\bar{a}_i - \bar{a}_{i-1}a_{i-1}, \ g_{2,i}^2 = -\bar{a}_{i-1}\bar{a}_{i-2}$$ for $i = 0, \ldots , m-1$.
We remark that the set $\{g^2_{r,i} \mid i=0, \ldots , m-1, r=0, 1, 2\}$ is a minimal set of uniform relations which generate  $I$.
Then, from \cite[Propositions 4.1, 5.1, 5.6, 6.2 and Theorem 7.2]{ST}, for all $m\pgq 1$, there is an element $\pi$ in $\HH^2(\Lambda)$ which is represented by the bimodule map $\pi : P^2 \to \Lambda$ in which the element $e_0\otimes e_0 \in \Lambda\mo(g^2_{1,0})\otimes\mt(g^2_{1,0})\Lambda$ has image $a_0\bar{a}_0 \in \Lambda$ and all other summands of $P^2$ have zero image. We label the idempotent generators of the summands of $P^2$ as follows: for each $i = 0, 1, \ldots , m-1$ and $r = 0, 1, 2$, write $e_i \otimes_{r,i} e_{i+2-2r}$ for the idempotent $\mo(g^2_{r,i})\otimes\mt(g^2_{r,i}) = e_i \otimes e_{i+2-2r}$ in the summand $\Lambda\mo(g^2_{r,i})\otimes\mt(g^2_{r,i})\Lambda$. When describing a map $P^2 \to \Lambda$, we omit summands whose image is zero.
Thus we may write $\pi$ as the map
$$\pi : e_0\otimes_{1,0} e_0 \mapsto a_0\bar{a}_0.$$
Now, $g^2_{1,0} = a_0\bar{a}_0 - \bar{a}_{m-1}a_{m-1}$. Since $\Lambda$ is Koszul, by \cite[Proposition 3.7]{BG}, the element $\pi$ gives rise to a unique formal deformation $\Lambda(T)$ of $\Lambda$, which, when we specialize the deformation parameter $T$ to $t \in K$, gives the algebra
$\Lambda(t) = K{\mathcal Q}/I(t)$, where $I(t)$ is the ideal of $K{\mathcal Q}$ generated by $a_ia_{i
+1}$, $\oa_{i-1}\oa_{i-2}$, $a_j\oa_j-\oa_{j-1}a_{j-1}$ and $(1-t)a_0\bar{a}_0 - \bar{a}_{m-1}a_{m-1}$ for $i = 0, \ldots , m-1, \ j = 1, \ldots , m-1$. We restrict ourselves to considering the case $t \neq 1$, since, if $t=1$, then the algebra $\Lambda(t)$ is not self-injective. In the case where $t = 0$, we recover the original algebra $\Lambda$. The algebra $\Lambda(t)$ for $t\in K\setminus\{1\}$ is a Koszul self-injective algebra, and we can easily verify that $\Lambda/\soc(\Lambda) \cong \Lambda(t)/\soc(\Lambda(t))$, so that $\Lambda(t)$ is a socle deformation of $\Lambda$.

\bigskip

This naturally leads us to introduce the algebra $\lq$ which we will study in this paper. Suppose $m \pgq 1$. For each ${\mathbf q} = (q_0, q_1, \ldots , q_{m-1}) \in (K^*)^m$, we define $\lq = K{\mathcal Q}/\iq$, where $\iq$ is the ideal of $K{\mathcal Q}$ generated by $$a_ia_{i+1}, \ \oa_{i-1}\oa_{i-2}, \ q_ia_i\oa_i-\oa_{i-1}a_{i-1} \mbox{ for } i = 0, \ldots , m-1.$$ Then $\Lambda(t) = \lq$ with ${\mathbf q} = (1-t, 1, \ldots , 1)$. We are assuming each $q_i$ is non-zero since we wish to study self-injective algebras. The algebra $\lq$ is a Koszul self-injective socle deformation of $\Lambda$, and $\lq = \Lambda$ when ${\mathbf q} = (1, 1, \ldots , 1)$.

Now, for $m \pgq 2$ and ${\mathbf q} = (1, q_1,1, \ldots , 1)$, the algebra $\Lambda_{(1, q_1,1, \ldots , 1)}$ comes from a formal deformation of $\Lambda$ via the element of $\HH^2(\Lambda)$ which is represented by the map $$\eta_1: P^2 \to \Lambda, e_1\otimes_{1,1} e_1 \mapsto a_1\bar{a}_1.$$ It can be easily verified using \cite{ST} that $\eta_1$ and $  \pi$ represent the same element in  $\HH^2(\Lambda)$. More generally, for $j = 1, \ldots , m-1$ the map
$$\eta_j: P^2 \to \Lambda, e_j\otimes_{1,j} e_j \mapsto a_j\bar{a}_j$$
also represents the element $\pi \in \HH^2(\Lambda)$. Thus the algebra $\Lambda_{(q_0, q_1, \ldots , q_{m-1})}$ comes from a formal deformation of $\Lambda$ by a scalar multiple of the element $\pi$. But $\Lambda_{(q_0, q_1, \ldots , q_{m-1})}$ can also be obtained from a  formal deformation of $\Lambda$ where we only replace the relation $g^2_{1,0} = a_0\bar{a}_0 - \bar{a}_{m-1}a_{m-1}$ by $(q_0\cdots q_{m-1})a_0\bar{a}_0 - \bar{a}_{m-1}a_{m-1}$ with $q_0\cdots q_{m-1} \in K^*$. Indeed we can give an explicit isomorphism $\Lambda_{(q_0, q_1, \ldots , q_{m-1})} \to \Lambda_{(q_0q_1\cdots q_{m-1}, 1, \ldots , 1)}$ as the algebra isomorphism induced by $a_i \mapsto q_0q_1\cdots q_ia_i, \oa_i \mapsto \oa_i$. Set $\zeta = q_0q_1\cdots q_{m-1} \in K^*$. Then $\Lambda_{(q_0, q_1, \ldots , q_{m-1})} \cong \Lambda_{(\zeta, 1, \ldots , 1)} = \Lambda(1-\zeta).$

\bigskip

However, there are other elements of $\HH^2(\Lambda)$ which we need to consider to see if they too give rise to a socle deformation of $\Lambda$. For $m \pgq 4$ and using \cite[Propositions 2.3, 2.4]{ST}, we have
$$\dim\HH^2(\Lambda) = \begin{cases}
1 & \mbox{if $m$ odd and $\car K \neq 2$}\\
2 & \mbox{if $m$ even, or if $m$ odd and $\car K = 2$.}
\end{cases}$$
Then, for $m \pgq 4$ and from \cite[Propositions 4.1, 5.1, 5.6]{ST}, $\HH^2(\Lambda)$ has basis
$$\begin{cases}
\{\pi\} & \mbox{if $m$ odd and $\car K \neq 2$}\\
\{\chi, \pi\} & \mbox{if $m$ even, or if $m$ odd and $\car K = 2$}
\end{cases}$$
where
$$\chi : e_i\otimes_{1,i} e_i \mapsto (-1)^ie_i \ \ \mbox{for $i = 0, \ldots , m-1.$}
$$

Let $m \pgq 4$ and let $\eta \in \HH^2(\Lambda)$. Then, by \cite[Proposition 3.7]{BG}, $\eta$ is the infinitesimal of the formal deformation of $\Lambda$ which, when the deformation parameter is specialized to $t \in K$, gives the algebra $ K{\mathcal Q}/J_{\eta}$, where $J_{\eta}$ is the ideal in $K{\mathcal Q}$ generated by $a_ia_{i+1}$, $\oa_{i-1}\oa_{i-2}$, $a_i\oa_i-\oa_{i-1}a_{i-1} - t\eta(e_i\otimes_{1,i} e_i)$ for $i = 0, \ldots , m-1$.

\bt
Let $m \pgq 4$, $\eta \in \HH^2(\Lambda)$ and $J_\eta$ as above. Then $ K{\mathcal Q}/J_{\eta}$ is a socle deformation of $\Lambda$ if and only if $\eta \in \sp\{\pi\}$.
\et

\begin{proof} As we have seen above, we can set $\eta=b_1\pi+b_2\chi$ for some $b_1, b_2$ in $K.$ Then the ideal $J_\eta$ is generated by $a_ia_{i+1}$, $\oa_{i-1}\oa_{i-2}$, $a_j\oa_j-\oa_{j-1}a_{j-1} - t(-1)^jb_2e_j$, $a_0\oa_0-\oa_{m-1}a_{m-1} - tb_2e_0-tb_1a_0\oa_0$ for $i = 0, \ldots , m-1$ and $j=1,\ldots,m-1$. Therefore the algebra $\tilde{\Lambda}:=K{\mathcal Q}/J_{\eta}$ has a $K$-basis given by $\set{e_i,a_i,\oa_i,a_i\oa_i;\ i=0,\ldots,m-1}.$

We first assume that $b_2\neq 0.$ Note that, for all $i=0,\ldots,m-2$, we have
\[ a_i\oa_ia_i=a_i\left(a_{i+1}\oa_{i+1}+(-1)^{i+1}tb_2e_{i+1}\right)=(-1)^{i+1}tb_2a_i \] and similarly, for $i=m-1,$  $a_{m-1}\oa_{m-1}a_{m-1}=(-1)^{m}tb_2a_{m-1}.$ Therefore, for any $i=0,\ldots,m-1,$ $\tilde{\Lambda}a_i$ has a $K$-basis given by $\set{a_i,\oa_i a_i}.$ Hence $\tilde{\Lambda}a_i$ is $2$-dimensional and it is easy to check that it is simple. We now show that the modules $\tilde{\Lambda}a_i$ for $i=0,\ldots,m-1$ are pairwise non-isomorphic. Suppose that there is a non-zero $\tilde{\Lambda}$-module morphism $f:\tilde{\Lambda}a_i\rightarrow \tilde{\Lambda}a_j.$ Then $e_if(a_i)=f(a_i)$ and $f(a_i)\in\sp\set{a_j,\oa_ja_j}$ so that $i=j$ or $i=j+1.$ If, moreover, $f$ is an isomorphism, we have a non-zero morphism $\tilde{\Lambda}a_j\rightarrow \tilde{\Lambda}a_i,$ and we get that $j=i$ or $j=i+1.$ Therefore there is an isomorphism $f:\tilde{\Lambda}a_i\rightarrow \tilde{\Lambda}a_j$ if and only if $i=j.$ Thus we have $m$ pairwise non-isomorphic $2$-dimensional simple $\tilde{\Lambda}$-modules, so that $\dim_K\soc(\tilde{\Lambda})\pgq 2m>m=\dim_K\soc(\Lambda)$. Hence $\tilde{\Lambda}/\soc(\tilde{\Lambda})$ is not isomorphic to $\Lambda/\soc(\Lambda)$ so that $\tilde{\Lambda}$ and $\Lambda$ are not socle equivalent.

Now assume that $b_2=0$ and $b_1\neq 0.$ Since a socle deformation of $\Lambda$ must be a self-injective algebra, necessarily $tb_1\neq 1.$ It is easy to check that the socle of $\tilde{\Lambda}$ is the submodule generated by the $a_i\oa_i$ for $i=0,\ldots,m-1$ so that $\tilde{\Lambda}/\soc(\tilde{\Lambda})\cong \Lambda/\soc(\Lambda)$, that is, $\tilde{\Lambda}$ is a socle deformation of $\Lambda.$
\end{proof}

Thus, for $m \pgq 4$, the socle deformations of $\Lambda$ which arise from formal deformations are precisely the algebras $\lq$, and the infinitesimal of the formal deformation is (a scalar multiple of) $\pi \in \HH^2(\Lambda)$.

\bigskip

For $m = 1, 2, 3$, there may be other socle deformations of $\lq$ which come from formal deformations. However, for $m=3$, it can be shown that there are no additional socle deformations arising in this way. But, for $m=2$, the elements $\pi_{2,-1}$ and $\pi_{2,1}$ in $\HH^2(\Lambda)$, which are given in \cite[Proposition 6.2]{ST} by
$$\begin{array}{l}
\pi_{2,-1} : e_0\otimes_{2,0} e_0 \mapsto a_0\oa_0\\
\pi_{2,1} : e_0\otimes_{0,0} e_0 \mapsto a_0\oa_0,
\end{array}$$
both give rise to the same socle deformation $\Lambda'$ of $\Lambda$, and, moreover, $\Lambda'$ is not isomorphic to $\lq$. We do not consider any additional socle deformations for $m = 1, 2$ in this paper.

\bigskip

Throughout this paper we suppose $m \pgq 1$, and consider the socle deformation $\lq$ of $\Lambda$. We write $\lq = K{\mathcal Q}/\iq$, where $\iq$ is the ideal generated by $a_ia_{i+1}$, $\oa_{i-1}\oa_{i-2}$ and  $q_ia_i\oa_i-\oa_{i-1}a_{i-1}$ for $i = 0, \ldots , m-1$, and
where ${\mathbf q} = (q_0, q_1, \ldots , q_{m-1}) \in (K^*)^m$ and $\zeta = q_0q_1\cdots q_{m-1} \in K^*$. In the case $m=1$, where the algebra $\Lambda_{(q)}$ was considered in \cite{BGMS}, different phenomena were exhibited depending on whether or not $q$ was a root of unity. Correspondingly, we will see in this current paper that we obtain different results depending on whether or not $\zeta$ is a root of unity.

\bigskip

\section{The graded centre of the Ext algebra of $\lq$}

\bigskip

We start by describing the Ext algebra $E(\lq)$. In Proposition \ref{prop:eltsincentre} we give some specific elements which lie in $Z_\gr(E(\lq))$. The remaining results lead to Theorem \ref{thm:structure} in which we prove that these elements generate the graded centre of the Ext algebra, thus enabling us to give a complete description of $Z_\gr(E(\lq))$. The algebras $\lq$ were studied in \cite{BGMS} in the case $m=1$; this
case is also included here.

The algebra $\lq$ is Koszul so, from \cite[Theorem 2.2]{GMV}, the Ext algebra $E(\lq)$
is the Koszul dual of $\lq$ and is given explicitly by quiver and relations
as $E(\lq) \cong K{\mathcal Q}^{\op}/\iq^\perp$, where ${\mathcal Q}$ is the quiver of $\lq$ and $\iq^\perp$ is the ideal of $K{\mathcal Q}^{\op}$ generated by
the orthogonal relations to those of $\iq$. Since left $K{\mathcal Q}^{\op}$-modules are right $K{\mathcal Q}$-modules, we may consider $E(\lq)$ as the quotient of  $K{\mathcal Q}$ by the ideal generated by $q_i^{-1}a_i\oa_i + \oa_{i-1}a_{i-1}$ for $i = 0, \ldots , m-1$, where we continue to write our paths from left to right. Let $\gamma^n_i$ denote the path $a_ia_{i+1}\cdots a_{i+n-1}$ of length $n$ in $K{\mathcal Q}$ which starts at the vertex $i$ and in which the subscripts are taken modulo $m$. Similarly we let $\delta^n_i$ denote the path $\oa_{i+n-1}\cdots \oa_{i+1}\oa_i$ of length $n$ in $K{\mathcal Q}$ which ends at the vertex $i$ and in which the subscripts are again taken modulo $m$. Thus a typical monomial in $E(\lq)$ has the form $\gamma_i^s\delta_j^t$ for some integers $s, t \pgq 0$ and $0 \ppq i, j \ppq m-1$. The algebra $E(\lq)$ is naturally graded by the length of paths. Note  that there is also a $\mathbb{Z}$-grading on the algebra $E(\lq)$ for which the degree of $a_i$ is $1$, the degree of $\bar{a}_i$ is $-1$ and the degree of $e_i$ is $0$ for all $i = 0, \ldots , m-1.$ Thus the length of a monomial $\gamma_i^s\delta_j^t$ is $s+t$ and its degree is $s-t.$ In particular, two such monomials with the same length have different degrees.

Let $z \in Z_\gr(E(\lq))$. Then $z = \sum_{i=0}^{m-1}e_ize_i$. For each $i = 0, \ldots , m-1$, a typical monomial in $e_iE(\lq)e_i$ has the form $\gamma_i^s\delta_i^t$ for some integers $s, t \pgq 0$ where $s \equiv t \pmod{m}$. Moreover, using the gradings mentioned above on $E(\lq)$, $Z_\gr(E(\lq))$ is generated by  elements which are both length homogeneous and degree homogeneous. So,  since we can assume that $z$ is homogeneous for both gradings,  we may write \[z = \sum_{i=0}^{m-1}c_i\gamma_i^{s_i}\delta_i^{t_i}\] where $c_i \in K$, $s_i, t_i \pgq 0$, $s_i \equiv t_i \pmod{m}$,  $s_i+t_i = s_0+t_0$ and $s_i-t_i=s_0-t_0$ for $i = 0, \ldots , m-1$. Keeping our convention on subscripts modulo $m$, we have $c_0=c_m$, $s_0 = s_m$ and $t_0 = t_m$.

\bigskip

Now, for $0 \ppq j \ppq m-1$, we have $$a_jz = c_{j+1}a_j\gamma_{j+1}^{s_{j+1}}\delta_{j+1}^{t_{j+1}} = c_{j+1}\gamma_j^{s_{j+1}+1}\delta_{j+1}^{t_{j+1}}$$
and
$$za_j = c_j\gamma_j^{s_j}\delta_j^{t_j}a_j = (-1)^{t_j}c_j(q_{j+1}\cdots q_{j+t_j})^{-1}\gamma_j^{s_{j}+1}\delta_{j+1}^{t_{j}}.$$
Since $a_jz = (-1)^{s_0+t_0}za_j$, we have, for all $j = 0, \ldots , m-1$, that either $c_j = 0$ or $c_{j+1} = (-1)^{s_j}c_j(q_{j+1}\cdots q_{j+t_j})^{-1}, s_j = s_{j+1}, t_j = t_{j+1}$.
If $c_j = 0$ for all $j = 0, \ldots , m-1$, then $z = 0$. So we assume now that $z \neq 0$. Then, for $j = 0, \ldots , m-1$, we have $s_j = s_0, t_j = t_0$ and $c_{j+1} = (-1)^{s_0}c_j(q_{j+1}\cdots q_{j+t_0})^{-1}\neq 0$. Thus  $z = \sum_{i=0}^{m-1}c_i\gamma_i^{s_0}\delta_i^{t_0} \neq 0$.

For $0 \ppq j \ppq m-1$, we have $$\oa_jz = c_j\oa_j\gamma_j^{s_0}\delta_j^{t_0} =
(-1)^{s_0}c_j(q_{j+1}\cdots q_{j+s_0})^{-1}\gamma_{j+1}^{s_0}\delta_j^{t_0+1}$$
and
$$z\oa_j = c_{j+1}\gamma_{j+1}^{s_0}\delta_{j+1}^{t_0}\oa_j = c_{j+1}\gamma_{j+1}^{s_0}\delta_{j}^{t_0+1}.$$
Since $\oa_jz = (-1)^{s_0+t_0}z\oa_j$, we also have, for all $j = 0, \ldots , m-1$, that $c_{j+1} = (-1)^{t_0}c_j(q_{j+1}\cdots q_{j+s_0})^{-1}$.

Thus $z = \sum_{i=0}^{m-1}c_i\gamma_i^{s_0}\delta_i^{t_0}$ with $c_{j+1} = (-1)^{s_0}c_j(q_{j+1}\cdots q_{j+t_0})^{-1} = (-1)^{t_0}c_j(q_{j+1}\cdots q_{j+s_0})^{-1}$ for $j = 0, \ldots , m-1$, and $s_0 \equiv t_0 \pmod{m}$.

From the equations $c_{j+1} = (-1)^{s_0}c_j(q_{j+1}\cdots q_{j+t_0})^{-1}$ we have that $$c_0 = (-1)^{ms_0}(q_0\cdots q_{t_0-1})^{-1}(q_1\cdots q_{t_0})^{-1}(q_{m-1}\cdots q_{m-2+t_0})^{-1}c_0.$$
Since $c_0 \neq 0$ and $\zeta = q_0\cdots q_{m-1}$ we get
$1 = (-1)^{ms_0}\zeta^{-t_0}$ so
$$\zeta^{t_0} = (-1)^{ms_0}.$$
In a similar way, the equations $c_{j+1} = (-1)^{t_0}c_j(q_{j+1}\cdots q_{j+s_0})^{-1}$
imply that $$\zeta^{s_0} = (-1)^{mt_0}.$$

\bigskip

It now follows immediately that if $\zeta$ is not a root of unity then $s_0 = t_0 = 0$, and so $c_j = c_0$ for all $j$. Hence $z = c_01$ with $c_0 \in K$. This gives the following result.

\bigskip

\bprop\label{prop:notunity}
If $\zeta$ is not a root of unity then $Z_\gr(E(\lq)) = K$.
\eprop

We now assume that $\zeta$ is a root of unity, and let $d \pgq 1$ be minimal such that $\zeta^d = 1$. We use the equations in $c_j$ and $c_{j+1}$ to write each $c_i$ in terms of $c_0$ for $i = 0, \ldots , m-1$. Thus we summarize the information about $z \in Z_\gr(E(\lq))$ as follows. We have
$z = \sum_{i=0}^{m-1}c_i\gamma_i^{s_0}\delta_i^{t_0}$ with $c_i = (-1)^{is_0}\prod_{k=1}^i(q_k\cdots q_{k+t_0-1})^{-1}c_0 = (-1)^{it_0}\prod_{k=1}^i(q_{k}\cdots q_{k+s_0-1})^{-1}c_0$, $\zeta^{s_0} = (-1)^{mt_0}$, $\zeta^{t_0} = (-1)^{ms_0}$ and $s_0 \equiv t_0 \pmod{m}$.

The next step is to verify that specific elements do indeed lie in the graded centre of the Ext algebra. The proof is straightforward and is omitted.

\bprop\label{prop:eltsincentre} Suppose that $\zeta$ is a primitive $d$-th root of unity.
\begin{enumerate}
\item Suppose that $m$ is even or $\car K = 2$. \\ Let $x = \sum_{i=0}^{m-1}\gamma_i^{dm}$, $y = \sum_{i=0}^{m-1}\delta_i^{dm}$ and $w = \sum_{i=0}^{m-1}(-1)^{id}\prod_{k=1}^i(q_k\cdots q_{k+d-1})^{-1}\gamma_i^{d}\delta_i^{d}$. Then $x, y, w \in Z_\gr(E(\lq))$.

    Moreover $w^m = \varepsilon xy$ where $\varepsilon = (-1)^{md/2}\prod_{l=1}^{m-1}\prod_{k=1}^{ld}(q_k\cdots q_{k+d-1})^{-1}$.
\item Suppose that $m$ is odd and $\car K \neq 2$.\\ Let $x=
\begin{cases}
\sum_{i=0}^{m-1}\gamma_i^{dm}&\text{if $d$ is even}\\\sum_{i=0}^{m-1}\gamma_i^{2dm}&\text{if $d$ is odd}
\end{cases},
$ \ \ let $y=
\begin{cases}
\sum_{i=0}^{m-1}\delta_i^{dm}&\text{if $d$ is even}\\\sum_{i=0}^{m-1}\delta_i^{2dm}&\text{if $d$ is odd}
\end{cases}$, \ \ and\\ let $w=\sum_{i=0}^{m-1}(-1)^{2\sigma i}\prod_{k=1}^i(q_k\cdots q_{k+\sigma d-1})^{-1}\gamma_i^{\sigma d}\delta_i^{\sigma d}$ where $\sigma=
\begin{cases}
1&\text{ if $d\equiv 0\pmod{4}$}\\\frac{1}{2}&\text{ if $d\equiv 2\pmod{4}$}\\2&\text{ if $d$ is odd.}
\end{cases}
$ Then $x, y, w \in Z_\gr(E(\lq))$.

Moreover $
\begin{cases}
 w^m=\varepsilon xy &\text{ if $d\equiv 0\pmod{4}$ or $d$ is odd }\\ w^{2m}=\varepsilon xy&\text{ if $d\equiv 2\pmod{4}$}
\end{cases}$ \\where $\varepsilon=
\begin{cases}
\prod_{l=1}^{m-1}\prod_{k=1}^{ld}(q_k\cdots q_{k+d-1})^{-1}&\text{ if $d\equiv 0\pmod{4}$}\\
\prod_{l=1}^{2m-1}\prod_{k=1}^{ld/2}(q_k\cdots q_{k+d-1})^{-1}&\text{ if $d\equiv 2\pmod{4}$}\\
\prod_{l=1}^{m-1}\prod_{k=1}^{2ld}(q_k\cdots q_{k+2d-1})^{-1}&\text{ if $d$ is odd.}\end{cases}
$
\end{enumerate}
\eprop

\bigskip

The main result of this section is Theorem \ref{thm:structure}, which shows that Proposition \ref{prop:eltsincentre} contains precisely the information needed to fully describe the graded centre $Z_\gr(E(\lq))$. Propositions \ref{prop:meven} and \ref{prop:modd} show that, where $\zeta$ is a root of unity, $Z_\gr(E(\lq))$ is indeed generated by $1, x, y$ and $w$ as a $K$-algebra. The next result, Lemma \ref{lem:nomorerelations}, is required to show that the only relation between the generators of $Z_\gr(E(\lq))$ is the relation of the form $w^p = \varepsilon xy$ as given in Proposition \ref{prop:eltsincentre}.

\bl\label{lem:nomorerelations}
With the notation of Proposition \ref{prop:eltsincentre}, suppose that $Z_\gr(E(\lq))$ is generated as an algebra by the elements $1, x, y$ and $w$ with homogeneous relation $w^p = \varepsilon xy$, for appropriate $\varepsilon \in K^*$ and positive integer $p$. Then $$Z_\gr(E(\lq)) = K[x,y,w]/\langle w^p - \varepsilon xy\rangle.$$
\el

\begin{proof}
Using the length grading on $E(\lq)$, we know that $Z_\gr(E(\lq))$ is a homogeneous quotient of $K[x,y,w]/\langle w^p - \varepsilon xy\rangle$, where $\varepsilon, p$ are as given in Proposition \ref{prop:eltsincentre}. Now, the elements $x^iy^{n-i}$, for $i = 0, \ldots , n$, are linearly independent in $E(\lq)$. So any additional relation in $Z_\gr(E(\lq))$ must be homogeneous of the form
$$f_0(x,y) + f_1(x,y)w + \cdots + f_{p-1}(x,y)w^{p-1} = 0$$
where $f_i(x,y) \in K[x,y]$ and $\deg(f_0(x,y)) = \deg(f_1(x,y)w) = \cdots = \deg(f_{p-1}(x,y)w^{p-1})$. Thus $\deg f_0(x,y) = \deg f_1(x,y) + |w|$ and, since $|x| = |y|$, there are integers $r, n$ with $\deg f_0(x,y) = n|x|$ and $\deg f_1(x,y) = r|x|$.

In the case $m$ even or $\car K = 2$ with $m \pgq 2$, we have $|x| = |y| = md, |w| = 2d$, which gives $nmd = rmd + 2d$ so that $2 = (n-r)m$. Since $m \pgq 2$, this implies $m = 2$ and $r = n-1$. Then $p=2$ and $|x| = |y| = |w| = 2d$. We may choose $n$ minimal so that
$f_0(x,y) + f_1(x,y)w = 0$ with $\deg f_0(x,y) = 2nd, \deg f_1(x,y) = 2(n-1)d$. Write
$f_0(x,y) = \sum_{i=0}^nb_ix^iy^{n-i}$ and $f_1(x,y) = \sum_{i=0}^{n-1}\tilde{b}_ix^iy^{n-i-1}$ with $b_i, \tilde{b}_i \in K$. Then $f_0^2(x,y) = f_1^2(x,y)w^2 = \varepsilon f_1^2(x,y)xy$. Equating coefficients of $x^{2n}$ and $y^{2n}$ gives that $b_0 = 0 = b_n$. Thus $f_0(x,y) = g(x,y)xy = \varepsilon^{-1}g(x,y)w^2$ for some $g(x,y) \in K[x,y]$. Hence $\varepsilon^{-1}g(x,y)w^2 + f_1(x,y)w = 0$ so that $\varepsilon^{-1}g(x,y)w + f_1(x,y) = 0$ which contradicts the minimality of $n$.

Now, suppose that $m$ is odd with $m \pgq 3$, and $\car K \neq 2$. If $d$ is even, we have $|x| = |y| = md$. If $d \equiv 0 \pmod{4}$, then $|w| = 2d$, which gives $nmd = rmd + 2d$ so that $2 = (n-r)m$. Since $m \pgq 3$ this has no solution. If $d \equiv 2 \pmod{4}$ then $|w| = d$, which gives $nmd = rmd + d$ so that $1 = (n-r)m$, and, again, this has no solution. Finally, if $d$ is odd, then $|x| = |y| = 2md$ and $|w| = 4d$. Hence $2nmd = 2rmd + 4d$ so that $2 = (n-r)m$, which also has no solution.

Finally we consider the case where $m = 1$. If $\car K = 2$, or if $\car K \neq 2$ and $d\equiv 0\pmod{4}$ or $d$ is odd, then we have, from Proposition \ref{prop:eltsincentre}, that $w = xy$, so that $Z_\gr(E(\lq)) = K[x,y]$. So suppose that $\car K \neq 2$ and $d\equiv 2\pmod{4}$. Then $|x| = |y| = d = |w|$ with $w^2 = \varepsilon xy$ where $\varepsilon = \prod_{k=1}^{d/2}(q_k\cdots q_{k+d-1})^{-1}$.
Then we have $nd = rd + d$ so that $r = n-1$. We may choose $n$ minimal so that
$f_0(x,y) + f_1(x,y)w = 0$ with $\deg f_0(x,y) = nd, \deg f_1(x,y) = (n-1)d$. Write
$f_0(x,y) = \sum_{i=0}^nb_ix^iy^{n-i}$ and $f_1(x,y) = \sum_{i=0}^{n-1}\tilde{b}_ix^iy^{n-i-1}$ with $b_i, \tilde{b}_i \in K$. We now apply the same argument as that used above for the case $m = 2$, to get a contradiction to the minimality of $n$.

\sloppy Thus there are no additional relations among the generators $x, y, w$, so it follows that $Z_\gr(E(\lq)) = K[x,y,w]/\langle w^p - \varepsilon xy\rangle.$
\end{proof}

The next stage is to determine $Z_\gr(E(\lq))$ in the case where $m$ is even or $\car K = 2$.

\bprop\label{prop:meven}
Suppose that $\zeta$ is a primitive $d$-th root of unity and that $m$ is even or $\car K = 2$. Then, keeping the notation of Proposition \ref{prop:eltsincentre}, $$Z_\gr(E(\lq)) = K[x, y, w]/\langle w^m - \varepsilon xy\rangle$$ where $\varepsilon = (-1)^{md/2}\prod_{l=1}^{m-1}\prod_{k=1}^{ld}(q_k\cdots q_{k+d-1})^{-1}$.
\eprop

\begin{proof}
If $m$ is even or $\car K = 2$, then $\zeta^{s_0} = (-1)^{mt_0} = 1$ and $\zeta^{t_0} = (-1)^{ms_0} = 1$. Thus $d|s_0$ and $d|t_0$. We also have that $s_0 \equiv t_0 \pmod{m}$ so, $t_0 = s_0 + rm$ for some integer $r$. We know
$c_1 = (-1)^{s_0}(q_1\cdots q_{t_0})^{-1}c_0 = (-1)^{t_0}(q_1\cdots q_{s_0})^{-1}c_0.$
If $m$ is even then $s_0$ and $t_0$ have the same parity, so we have that $(-1)^{s_0} = (-1)^{t_0}$. Hence $q_1\cdots q_{t_0} = q_1\cdots q_{s_0}$. Thus, if $t_0 \pgq s_0$, we have $q_{s_0+1}\cdots q_{t_0} = 1$, and if $s_0 \pgq t_0$, then we have $q_{t_0+1}\cdots q_{s_0} = 1$. Hence, in both cases, we have $\zeta^r = 1$ and $d|r$. Thus $t_0 = s_0 + hdm$ for some integer $h$.

Write  $z \in Z_\gr(E(\lq))$ as $z = \sum_{i=0}^{m-1}c_i\gamma_i^{s_0}\delta_i^{t_0}$. Suppose first that $s_0 = t_0 = 0$. Then $c_i = c_0$ for $i = 1, \ldots , m-1$ so that $z = c_01$. Now suppose that $s_0 = 0$ but $t_0 \neq 0$. Then $z = \sum_{i=0}^{m-1}c_i\delta_i^{t_0}$ with $t_0 = hdm$ for some $h \geq 1$ and $c_i = (-1)^{it_0}c_0$ for $i = 0, \ldots , m-1$. Since $m$ is even or $\car K = 2$, we have $c_i = c_0$ for all $i$ and so $z = c_0\sum_{i=0}^{m-1}\delta_i^{hdm} = c_0(\sum_{i=0}^{m-1}\delta_i^{dm})^h = c_0y^h$.
Similarly, if $t_0 = 0$ but $s_0 \neq 0$, then $s_0 = hdm$ for some $h \geq 1$ and $z = c_0x^h$.

So suppose now that $s_0 \neq 0$ and $t_0 \neq 0$. Without loss of generality, assume that $t_0 \pgq s_0$ so $t_0 = s_0 + hdm$ for some integer $h \pgq 0$. Recalling that $(-1)^{s_0} = (-1)^{t_0}$, then
$$\begin{array}{rl}
z = & \sum_{i=0}^{m-1}(-1)^{it_0}\prod_{k=1}^i(q_{k}\cdots q_{k+s_0-1})^{-1}c_0\gamma_i^{s_0}\delta_i^{t_0} \\
= & \sum_{i=0}^{m-1}(-1)^{is_0}\prod_{k=1}^i(q_{k}\cdots q_{k+s_0-1})^{-1}c_0\gamma_i^{s_0}\delta_i^{s_0}\delta_i^{hdm}\\
= & c_0(\sum_{i=0}^{m-1}(-1)^{is_0}\prod_{k=1}^i(q_{k}\cdots q_{k+s_0-1})^{-1}\gamma_i^{s_0}\delta_i^{s_0})(\sum_{i=0}^{m-1}\delta_i^{dm})^h\\
= & c_0(\sum_{i=0}^{m-1}(-1)^{is_0}\prod_{k=1}^i(q_{k}\cdots q_{k+s_0-1})^{-1}\gamma_i^{s_0}\delta_i^{s_0})y^h.
\end{array}$$
Write $s_0 = \alpha dm + s$ with $0 \ppq s < dm$. Then (using again that $m$ is even or $\car K = 2$), we have $(-1)^{s_0} = (-1)^s$, and $q_k\cdots q_{k+s_0-1} = \zeta^{\alpha d}q_k\cdots q_{k+s-1} = q_k\cdots q_{k+s-1}$. Also, $\gamma_i^{s_0}\delta_i^{s_0} = \gamma_i^s(\gamma_{i+s}^{s_0 -s}\delta_{i+s}^{s_0-s})\delta_i^s = \gamma_i^s(\sum_{j=0}^{m-1}\gamma_{j}^{\alpha dm}\delta_{j}^{\alpha dm})\delta_i^s$. Now $\sum_{j=0}^{m-1}\gamma_{j}^{\alpha dm}\delta_{j}^{\alpha dm} = (\sum_{j=0}^{m-1}\gamma_{j}^{\alpha dm})(\sum_{j=0}^{m-1}\delta_{j}^{\alpha dm})
= x^{\alpha}y^{\alpha}$. So $\gamma_i^{s_0}\delta_i^{s_0} = \gamma_i^s\delta_i^sx^{\alpha}y^{\alpha}$ by Proposition \ref{prop:eltsincentre}. Thus it is sufficient to consider $$z = \sum_{i=0}^{m-1}(-1)^{is}\prod_{k=1}^i(q_{k}\cdots q_{k+s-1})^{-1}\gamma_i^{s}\delta_i^{s}$$ where $0 \ppq s \ppq dm-1$.

Now $d|s_0$ so $d|s$, and thus $s \in \{0, d, 2d, \ldots , (m-1)d\}$. Let $s = jd$ and define $$z_j = \sum_{i=0}^{m-1}(-1)^{ijd}\prod_{k=1}^i(q_k\cdots q_{k+jd-1})^{-1}\gamma_i^{jd}\delta_i^{jd}$$ for $j = 0, 1, \ldots , m$. If $j = 0$ then $z_0 = 1$; if $j = 1$ then $z_1 = w$, and if $j=m$ then $z_m = \sum_{i=0}^{m-1}(-1)^{imd}\prod_{k=1}^i(q_k\cdots q_{k+md-1})^{-1}\gamma_i^{md}\delta_i^{md} = \sum_{i=0}^{m-1}\gamma_i^{md}\delta_i^{md} = xy$. Moreover, it is easy to verify that $z_jw = (-1)^{jd}\prod_{k=1}^{jd}(q_k \cdots q_{k+d-1})^{-1}z_{j+1}$ for $j = 0, 1, \ldots , m-1$. We also have that $w^j = (-1)^{\sum_{i=1}^{j-1}id}(\prod_{l=1}^{j-1}\prod_{k=1}^{ld}(q_k\cdots q_{k+d-1})^{-1})z_j$ for $j = 0, 1, \ldots , m-1$. Hence $Z_\gr(E(\lq))$ is generated as an algebra by $1, x, y, w$ with $w^m = (-1)^{md/2}(\prod_{l=1}^{m-1}\prod_{k=1}^{ld}(q_k\cdots q_{k+d-1})^{-1})xy = \varepsilon xy$. The result now follows from Lemma \ref{lem:nomorerelations}.
\end{proof}

We now consider the case where $m$ is odd and $\car K \neq 2$.

\bprop\label{prop:modd}
Suppose that $\zeta$ is a primitive $d$-th root of unity, that $m$ is odd and $\car K \neq 2$. Then, keeping the notation of Proposition \ref{prop:eltsincentre}, $$Z_\gr(E(\lq)) =
\begin{cases}
K[x, y, w]/\langle w^m - \varepsilon xy \rangle & \text{ if $d\equiv 0\pmod{4}$ or $d$ is odd }\\
K[x, y, w]/\langle w^{2m} - \varepsilon xy \rangle & \text{ if $d\equiv 2\pmod{4}$}
\end{cases}$$
where $\varepsilon=
\begin{cases}
\prod_{l=1}^{m-1}\prod_{k=1}^{ld}(q_k\cdots q_{k+d-1})^{-1}&\text{ if $d\equiv 0\pmod{4}$}\\
\prod_{l=1}^{2m-1}\prod_{k=1}^{ld/2}(q_k\cdots q_{k+d-1})^{-1}&\text{ if $d\equiv 2\pmod{4}$}\\
\prod_{l=1}^{m-1}\prod_{k=1}^{2ld}(q_k\cdots q_{k+2d-1})^{-1}&\text{ if $d$ is odd.}\end{cases}$
\eprop

\begin{proof}
From the conditions $\zeta^{s_0} = (-1)^{mt_0}$ and $\zeta^{t_0} = (-1)^{ms_0}$, we get $\zeta^{2s_0} = 1 = \zeta^{2t_0}$ which gives $d|2s_0$ and $d|2t_0$.
We also have that $s_0 \equiv t_0 \pmod{m}$ so, $t_0 = s_0 + rm$ for some integer $r$. We know $c_1 = (-1)^{s_0}(q_1\cdots q_{t_0})^{-1}c_0 = (-1)^{t_0}(q_1\cdots q_{s_0})^{-1}c_0.$ If $t_0 \pgq s_0$, we have $q_{s_0+1}\cdots q_{t_0} = (-1)^{t_0-s_0}$, and, if $s_0 \pgq t_0$, we have $q_{t_0+1}\cdots q_{s_0} =(-1)^{s_0-t_0} = (-1)^{t_0-s_0}$. Hence, in both cases, $\zeta^r = (-1)^{t_0-s_0}$. Thus $d|2r$ and so $dm|2(t_0-s_0)$.
Now write $s_0 = \alpha dm + s$ and $t_0 = \beta dm + t$ with $0 \ppq s, t < dm$. Then $dm|2(t-s)$. Moreover, we may assume without loss of generality, that $t \pgq s$ so that $2(t-s)\in\set{0,dm}$.

We wish to show that $2(t-s) = 0$. So, we assume first that $2(t-s) = dm$ and aim for a contradiction. Since $m$ is odd, $2(t-s) = dm$ implies that $d$ is even. In particular, $t_0-s_0$ and $t-s$ have the same parity.  Moreover, $(-1)^{t_0-s_0} = \zeta^r = \zeta^{(t_0-s_0)/m} = \zeta^{(t-s)/m} = \zeta^{d/2} = -1$. Thus $t-s$ is odd and $\frac{d}{2}$ is odd. But $m$ is odd, so we can also use our first conditions to get $(-1)^{s_0+t_0} = (-1)^{m(s_0+t_0)} = (-1)^{ms_0}(-1)^{mt_0} = \zeta^{s_0+t_0} = \zeta^{s+t} = \zeta^{2s+(t-s)} = \zeta^{2s+(dm/2)} = \zeta^{2s}(-1)^m = -\zeta^{2s}$. Thus, squaring this identity gives $1 = \zeta^{4s}$ so that $d|4s$ and hence $\frac{d}{2}|2s$. But $\frac{d}{2}$ is odd so $\frac{d}{2}|s$ and we may set $s=\frac{d}{2}l$ for some integer $l$. However, if $s$ and therefore $l$ are both even, we get $1=(-1)^s=(-1)^{s_0}=\zeta^{t_0}=\zeta^t=\zeta^{(l+m)d/2}=(-1)^{l+m}=-1$, a contradiction, and if $s$ and therefore $l$ are both odd, then $t$ is even and we get $1=(-1)^t=(-1)^{t_0}=\zeta^{s_0}=\zeta^s=\zeta^{ld/2}=(-1)^l=-1,$ a contradiction. Thus $2(t-s) \neq dm$.

Therefore $2(t-s)=0$ and hence $t=s$. In this case, $1 = \zeta^{(\alpha-\beta)dm} = \zeta^{s_0-t_0}=(-1)^{m(t_0-s_0)} = (-1)^{t_0-s_0}=(-1)^{(\beta-\alpha)dm}$ so $\alpha d$ and $\beta d$ have same parity. Moreover, $\zeta^s=\zeta^{s_0}=(-1)^{t_0}$ so $d|2s$ with $0\ppq 2s<2dm.$ Hence $2s=ld$ for some integer $l$ with $0\ppq l < 2m.$
If $2s=ld$ with $l$ odd, then $d$ is even and $-1=(-1)^l = (\zeta^{d/2})^l = \zeta^s = \zeta^{s_0} = (-1)^{t_0} = (-1)^s=(-1)^{ld/2}=(-1)^{d/2}$ so $\frac{d}{2}$ is odd. On the other hand, if $2s=ld$ with $l$ even, then $1 = \zeta^{ld/2} = \zeta^s = \zeta^{s_0} = (-1)^{t_0} = (-1)^{s+\beta d} = (-1)^{(l/2+\beta)d}$ so $(\beta+\frac{l}{2})d$ is even, and consequently $(\alpha+\frac{l}{2})d$ is even. In this case, we also have that $t_0$ is even.

We are now able to describe the elements of $Z_\gr(E(\lq))$. Recall that
a typical homogeneous non-zero element $z\in Z_\gr(E(\lq))$ has the form $z = \sum_{i=0}^{m-1}(-1)^{it_0}\prod_{k=1}^i(q_{k}\cdots q_{k+s_0-1})^{-1}c_0\gamma_i^{s_0}\delta_i^{t_0} = \sum_{i=0}^{m-1}(-1)^{it_0}\prod_{k=1}^i(q_{k}\cdots q_{k+s-1})^{-1}c_0\gamma_i^{s_0}\delta_i^{t_0}$ for some $c_0 \in K^*$. We keep the notation of Proposition \ref{prop:eltsincentre} when referring to $x, y, w$.

If $d$ is odd, then $l$ and $t_0$ are even, so $z = \sum_{i=0}^{m-1}\prod_{k=1}^i(q_k\cdots q_{k+ld/2-1})^{-1}c_0\gamma_i^{\alpha dm+ld/2}\delta_i^{\beta dm+ld/2}$ with $\alpha, \beta$ integers such that $\alpha+\frac{l}{2}$ and $\beta+\frac{l}{2}$ are even and $0\ppq \frac{l}{2}\ppq m-1.$ If $\alpha$ is even and we let $\frac{l}{2}= 2L$, then $z$ is a scalar multiple of $x^{\alpha/2}y^{\beta/2}w^{L}$. If $\alpha$ is odd and we let $\frac{l}{2}= L$ then $z$ is a scalar multiple of $x^{(\alpha-1)/2}y^{(\beta-1)/2}w^{(m+L)/2}$.

If $d$ is even with $d\equiv 0\pmod{4}$, then $l$ and $t_0$ are even. Then $0\ppq \frac{l}{2}\ppq m-1$ and $z = \sum_{i=0}^{m-1}\prod_{k=1}^i(q_k\cdots q_{k+ld/2-1})^{-1}c_0\gamma_i^{\alpha dm+ld/2}\delta_i^{\beta dm+ld/2}$ with $\alpha$, $\beta$ integers. Hence $z$ is a scalar multiple of $x^\alpha y^\beta w^{l/2}$.

Finally, if $d$ is even with $d\equiv 2\pmod{4}$, then $l$ and $t_0$ have the same parity, so that $z = \sum_{i=0}^{m-1}(-1)^{li}\prod_{k=1}^i(q_k\cdots q_{k+l d/2-1})^{-1}c_0\gamma_i^{\alpha dm+ld/2}\delta_i^{\beta dm+ld/2}$ with $\alpha,$ $\beta$ and $l$ integers such that $0\ppq l\ppq 2m-1.$ In this case, $z$ is a scalar multiple of $x^\alpha y^\beta w^l.$

Thus $Z_\gr(E(\lq))$ is generated as an algebra by $1, x, y$ and $w$, where $x, y, w$ are as in Proposition \ref{prop:eltsincentre}. It remains to verify the relations of the form $w^p = \varepsilon xy$, for appropriate $\varepsilon \in K^*$ and positive integer $p$. The proofs are straightforward and left to the reader. The final description now follows from Lemma \ref{lem:nomorerelations}.
\end{proof}

We summarize Propositions \ref{prop:notunity}, \ref{prop:meven} and \ref{prop:modd} in the following result.

\bt\label{thm:structure}
Let ${\mathbf q} = (q_0, q_1, \ldots , q_{m-1}) \in (K^*)^m$ and let $\zeta = q_0q_1\cdots q_{m-1}$. If $\zeta$ is not a root of unity then $Z_\gr(E(\lq)) = K$.
Now suppose that $\zeta$ is a primitive $d$-th root of unity.
\begin{enumerate}
\item If $m$ is even or if $\car K = 2$, then $$Z_\gr(E(\lq)) = K[x, y, w]/\langle w^m - \varepsilon xy\rangle,$$
     where $\varepsilon = (-1)^{md/2}\prod_{l=1}^{m-1}\prod_{k=1}^{ld}(q_k\cdots q_{k+d-1})^{-1}$.
\item If $m$ is odd and $\car K \neq 2$, then
$$Z_\gr(E(\lq)) =
\begin{cases}
K[x, y, w]/\langle w^m - \varepsilon xy \rangle & \text{ if $d\equiv 0\pmod{4}$ or $d$ is odd }\\
K[x, y, w]/\langle w^{2m} - \varepsilon xy \rangle & \text{ if $d\equiv 2\pmod{4}$}
\end{cases}$$
where $\varepsilon=
\begin{cases}
\prod_{l=1}^{m-1}\prod_{k=1}^{ld}(q_k\cdots q_{k+d-1})^{-1}&\text{ if $d\equiv 0\pmod{4}$}\\
\prod_{l=1}^{2m-1}\prod_{k=1}^{ld/2}(q_k\cdots q_{k+d-1})^{-1}&\text{ if $d\equiv 2\pmod{4}$}\\
\prod_{l=1}^{m-1}\prod_{k=1}^{2ld}(q_k\cdots q_{k+2d-1})^{-1}&\text{ if $d$ is odd.}\end{cases}$
\end{enumerate}
\et

\bigskip

\section{The Hochschild cohomology ring modulo nilpotence of $\lq$}

\bigskip

We begin with the following corollary of Theorem \ref{thm:structure}.

\bc
Let ${\mathbf q} = (q_0, q_1, \ldots , q_{m-1}) \in (K^*)^m$ and let $\zeta = q_0q_1\cdots q_{m-1}$. Then $E(\lq)$ is finitely generated over $Z_\gr(E(\lq))$ if and only if $\zeta$ is a root of unity.
\ec

\begin{proof}
Since $\lq$ is a Koszul algebra, $E(\lq)$ is generated as a $K$-algebra in degrees 0 and 1. If $\zeta$ is not a root of unity, then $E(\lq)$ is not a finitely generated module over $Z_\gr(E(\lq))$ since $E(\lq)$ is an infinite-dimensional vector space. If $\zeta$ is a root of unity, then it is straightforward to verify that the set $\{\gamma_i^s\delta_j^t \mid 0 \ppq i, j \ppq m-1, 0 \ppq s, t \ppq |x|\}$
is a sufficient (but not necessarily minimal) generating set for $E(\lq)$ as a $Z_\gr(E(\lq))$-module, where the degree of $x$ is as given in Proposition \ref{prop:eltsincentre}.
\end{proof}

Using \cite{BGSS,SS}, $\HH^*(\Lambda_q)/\N \cong Z_\gr(E(\lq))/\N_Z$, where $\N_Z$ denotes the ideal of $Z_\gr(E(\lq))$ which is generated by all nilpotent elements. It is clear from Theorem \ref{thm:structure} that $\N_Z = 0$ so that $\HH^*(\Lambda_q)/\N \cong Z_\gr(E(\lq))$. Thus we have the following result.

\bt\label{thm:HH}
Let ${\mathbf q} = (q_0, q_1, \ldots , q_{m-1}) \in (K^*)^m$ and let $\zeta = q_0q_1\cdots q_{m-1}$. If $\zeta$ is not a root of unity then $\HH^*(\lq)/\N \cong K$.
If $\zeta$ is a root of unity, then $\HH^*(\lq)/\N$ is a finitely generated commutative $K$-algebra of Krull dimension 2.
\et

In particular, the conjecture of \cite{SS} holds for the class of algebras $\lq$ for all ${\mathbf q} \in (K^*)^m$ .

\end{document}